\documentclass[12pt,leqno]{article}
\usepackage{amsmath, amsthm, amsfonts, amssymb, color}
\usepackage{mathrsfs}
\theoremstyle{definition}

\newcommand{\scr}[1]{\mathscr #1}
\definecolor{wco}{rgb}{0.5,0.2,0.3}

\numberwithin{equation}{section} \theoremstyle{remark}

\newcommand{\ua}{\uparrow}

\title{{\bf Equivalent Semigroup Properties for   Curvature-Dimension Condition}\footnote{Supported in
 part by WIMCS and SRFDP.}
}
\author{
{\bf Feng-Yu Wang}\\
\footnotesize{School of Mathematical Sciences,
Beijing Normal
University, Beijing 100875, China}\\
\footnotesize{and}\\ \footnotesize{Department of Mathematics,
Swansea University, Singleton Park, SA2 8PP, UK}\\ \footnotesize{Email: wangfy@bnu.edu.cn;
F.Y.Wang@swansea.ac.uk}}
\begin{document}
\def\R{\mathbb R}  \def\ff{\frac} \def\ss{\sqrt} \def\BB{\mathbb
B}
\def\N{\mathbb N} \def\kk{\kappa} \def\m{{\bf m}}
\def\dd{\delta} \def\DD{\Delta} \def\vv{\varepsilon} \def\rr{\rho}
\def\<{\langle} \def\>{\rangle} \def\GG{\Gamma} \def\gg{\gamma}
  \def\nn{\nabla} \def\pp{\partial} \def\tt{\tilde}
\def\d{\text{\rm{d}}} \def\bb{\beta} \def\aa{\alpha} \def\D{\scr D}
\def\E{\mathbb E} \def\si{\sigma} \def\ess{\text{\rm{ess}}}
\def\beg{\begin} \def\beq{\begin{equation}}  \def\F{\scr F}
\def\Ric{\text{\rm{Ric}}} \def\Hess{\text{\rm{Hess}}}\def\B{\mathbb B}
\def\e{\text{\rm{e}}} \def\ua{\underline a} \def\OO{\Omega} \def\sE{\scr E}
\def\oo{\omega}     \def\tt{\tilde} \def\Ric{\text{\rm{Ric}}}
\def\cut{\text{\rm{cut}}} \def\P{\mathbb P} \def\ifn{I_n(f^{\bigotimes n})}
\def\C{\scr C}      \def\aaa{\mathbf{r}}     \def\r{r}
\def\gap{\text{\rm{gap}}} \def\prr{\pi_{{\bf m},\varrho}}  \def\r{\mathbf r}
\def\Z{\mathbb Z} \def\vrr{\varrho} \def\ll{\lambda}
\def\L{\scr L}\def\Tt{\tt} \def\TT{\tt}\def\II{\mathbb I}
\def\i{{\rm i}}\def\Sect{{\rm Sect}}

\maketitle
\begin{abstract}  Some equivalent gradient and Harnack inequalities of a diffusion semigroup are presented for the curvature-dimension condition of the associated generator. As applications, the first eigenvalue,  the log-Harnack inequality, the heat kernel estimates, and the HWI inequality are derived by using the curvature-dimension condition. The transportation inequality for diffusion semigroups is also investigated.
\end{abstract} \noindent
 \vskip 2cm

\section{Introduction}
Let $M$ be  a $d$-dimensional complete connected Riemannian manifold  without boundary or with a convex boundary $\pp M$. Let $P_t$ be the (Neumann if $\pp M\ne\emptyset$) semigroup generated by $L=\DD+Z$ for a $C^1$-vector field $Z$ on $M$. To describe  analytic properties of $P_t$,  the following curvature-dimension condition of Bakry-Emery \cite{BE} plays a very important role:

\beq\label{CD} \ff 12 L |\nn f|^2 -\<\nn Lf,\nn f\>\ge -K |\nn f|^2 +\ff 1 n (Lf)^2,\ \ f\in C^\infty(M),\end{equation} where $-K\in \R$ and $n\ge d$ provide   a curvature lower bound and   a dimension upper bound of $L$ respectively. When $Z=0$ this condition is equivalent to $\Ric\ge -K$, where $\Ric$ is the Ricci curvature. In this case (\ref{CD}) holds for $n=d$. When $Z\ne 0$, $n$ is essentially larger than $d$. Indeed,  (\ref{CD}) is equivalent to

\beq\label{CD'}\Ric(U,U)-\<\nn_U Z,U\>\ge -K|U|^2 - \ff{\<Z,U\>^2}{n-d},\ \ U\in TM.\end{equation}
In particular, when $n=\infty$,  (\ref{CD}) reduces to the curvature condition

\beq\label{C}  \Ric(U,U)-\<\nn_U \nn Z, U\>\ge -K |U|^2,\ \ U\in TM.\end{equation}

There are a number of equivalent semigroup inequalities for the curvature condition (\ref{C}), including  gradient inequalities, Poincar\'e/log-Sobolev inequalities, the dimension-free and logarithmic Harnack inequalities,  and Wasserstein (or transportation-cost) inequalities,  see e.g.
\cite{B, L, RS, W04, W10} and references within for details.

When $n<\infty$, the curvature-dimension condition (\ref{CD}) has been used in the study of the Sobolev inequality, the first eigenvalue and the diameter estimates, and Li-Yau type Harnack inequalities. Besides the above mentioned references, we refer to \cite{BQ99, BQ00, S} and references within for detailed applications of the curvature-dimension condition.
On the other hand, however, unlike for (\ref{C}), there is no any known equivalent semigroup inequalities for the curvature-dimension condition (\ref{CD}) with finite $n$.  The purpose of this note is to find  inequalities of $P_t$ which are equivalent to (\ref{CD}), and to make further applications of these equivalent inequalities.

Let $\D_0$ be the set of all smooth functions on $M$ with compact support and satisfying the Neumann boundary condition provided $\pp M\ne \emptyset$. Recall that throughout the paper $\pp M$ is assumed to be convex if it exists. Let $\rr$ be the Riemannian distance on $M$.

\beg{thm}\label{T1.1} Each of the following statements is equivalent to $(\ref{CD})$:

\beg{enumerate} \item[$(1)$] \ $|\nn P_t f|^2 \le \e^{2Kt} P_t |\nn f|^2-\ff 2 n   \int_0^t\e^{2Ks}P_s (P_{t-s}Lf)^2 \d s,\ \ f\in \D_0, t\ge 0.$ \item[$(2)$] \ $|\nn P_t f|^2 \le \e^{2Kt} P_t |\nn f|^2-\ff{\e^{2Kt}-1}{Kn} (P_t Lf)^2, \ \ f\in \D_0, t\ge 0.$
\item[$(3)$] \ $P_tf^2 -(P_t f)^2\le \ff{\e^{2Kt}-1}K P_t|\nn f|^2 -\ff{\e^{2Kt}-1-2Kt}{K^2 n} (P_t Lf)^2,\ \ f\in \D_0, t\ge 0.$
\item[$(4)$] \ $P_tf^2 -(P_t f)^2\ge \ff{1-\e^{-2Kt}}K |\nn P_tf|^2 +\ff{\e^{-2Kt}-1+2Kt}{K^2 n} (P_t Lf)^2,\ \ f\in \D_0, t\ge 0.$
\item[$(5)$] \ $\e^{Kt}P_t |\nn f|\ge |\nn P_t f| +\ff 1 {n-d} \int_0^t \e^{Ks} P_s \ff {\<Z, \nn P_{t-s} f\>^2}{|\nn P_{t-s} f|} \d s,\ \ f\in \D_0, t\ge 0.$
\item[$(6)$] For any $t>0$ and increasing $\varphi\in C^1([0,t])$ with $\varphi(0)=0$ and $\varphi'(0)=1$, the log-Harnack inequality
    $$P_{\varphi(t)} \log f(y)\le \log P_t f(x)
+ \ff{\rr(x,y)^2}{4\int_0^t \e^{-2K\varphi(s)} )\d s}
+\ff{Kn}4 \int_0^t \ff{(\varphi'(s)-1)^2\d s}{1-\e^{-2K\varphi(s)}}$$ holds for any positive function $f$ with $\inf f>0$ and all  $x,y\in M.$\end{enumerate}
\end{thm}
We remark that according to \cite[Theorem 1.2]{W10}, at least for compact manifolds and a class of non-compact manifolds, any of statements (1)-(6) implies that $\pp M$ is convex if exists. Therefore, our assumption on the boundary is essential.

 Now, we consider applications of the above equivalent inequalities.
We first present some consequences of   (6) for heat kernel bounds and HWI inequalities. According to Li-Yau's Harnack inequality \cite{LY, BQ99}, if (\ref{CD}) holds then   $P_t$ can be dominated by $P_{t+s}$ for $s,t>0.$ A nice point of (6)   is that we are also  able to dominate $P_{t+s}$ by $P_t$ with help of the logarithmic function. With concrete choices of $\varphi$ we have  the following explicit log-Harnack inequalities.

\beg{cor} \label{C1.3}  If $(\ref{CD})$ holds, then for any $s\ge 0, t>0,$

\beq\label{H1} P_{t+s}\log f(y)\le \log P_t f(x) +\ff{K(t+2s)\rr(x,y)^2}{2t(1-\e^{-2K(t+s)})} +\ff {n Ks^2}{2t(1-\e^{-Kt})}, \end{equation} and
\beq\label{H2} P_{t}\log f(y)\le \log P_{t+s} f(x) +\ff{K\rr(x,y)^2}{2(1-\e^{-2Kt})+4Ks\e^{-2Kt}} +\ff {Kns}{4(1-\e^{-2Kt})} \end{equation} hold for $x,y\in M$ and bounded measurable function $f$ with $\inf f>0$.
\end{cor}

As shown in the proof of \cite[Proposition 2.4(2)]{W10}, it is easy to see that for any $t>0, s\ge 0$ and $x,y\in M$, (\ref{H1}) and (\ref{H2}) are equivalent to the following heat kernel inequalities (\ref{H1'}) and (\ref{H2'}) respectively, where $\nu$ is a measure equivalent to $\d x$ and $p_t^\nu$ is the heat kernel of $P_t$ w.r.t. $\nu$:

\beq\label{H1'} \int_M p_{t+s}^\nu(y,z)\log \ff{p_{t+s}^\nu(y,z)}{p_t^\nu(x,z)}\,\nu(\d z)\le \ff{K(t+2s)\rr(x,y)^2}{2t(1-\e^{-2K(t+s)})} +\ff {n Ks^2}{2t(1-\e^{-Kt})},\end{equation}
\beq\label{H2'} \int_M p_{t}^\nu(y,z)\log \ff{p_{t}^\nu(y,z)}{p_{t+s}^\nu(x,z)}\,\nu(\d z)\le\ff{K\rr(x,y)^2}{2(1-\e^{-2Kt})+4Ks\e^{-2Kt}} +\ff {Kns}{4(1-\e^{-2Kt})}.\end{equation}
In particular, when $P_t$ is symmetric w.r.t a probability measure $\mu$,   we have the following  heat kernel lower bound.

\beg{cor}\label{C1.3'} Let $Z=\nn V$ such that $\mu(\d x):=\e^{V(x)}\d x$ is a probability measure, and let $p_t(x,y)$ be the heat kernel of $P_t$ w.r.t. $\mu$. Then $(\ref{C})$ and hence $(\ref{CD})$ implies
\beq\label{Heat} p_t(x,y)\ge  \exp\bigg[-\ff{K\rr(x,y)^2}{2(1-\e^{-Kt})}\bigg],\ \ x,y\in M, t>0.\end{equation}\end{cor}
We remark that (\ref{Heat}) is new. Known heat kernel lower bounds derived from Li-Yau's Harnack inequality are dimension-dependent, and decay to zero as the dimension goes to infinity provided $K>0$, see e.g. \cite[Corollary 3.9]{St} and \cite[(13)]{BQ99}.

Moreover, following the line of \cite{BGL}, we use the  log-Harnack inequality (\ref{H2}) to establish the HWI inequality. Again let $Z=\nn V$ such that
$\mu(\d x):=\e^{V(x)}\d x$ is a probability   measure. Recall that for any non-negative measurable function ${\bf c}$ on
$M\times M$, and for any $p\ge 1$,  the $L^p$-transportation cost induced by cost function ${\bf c}$ is

$$W_p^{\bf c}(\mu_1,\mu_2)= \inf_{\pi\in\scr C(\mu_1,\mu_2)} \pi({\bf c}^p)^{1/p},\ \ \mu_1,\mu_2\in\scr P(M),$$ where $\scr P(M)$ is the set of all probability measures on $M$ and $\scr C(\mu_1,\mu_2)$ is the set of all couplings for $\mu_1$ and $\mu_2$.

\beg{cor}\label{C1.4} Let $Z=\nn V$ such that $\mu(\d x):=\e^{V(x)}\d x$ is a probability measure. If $(\ref{CD})$ holds, then for any $f\in C^1(M)$  with $\mu(f^2)=1$,

\beg{equation}\label{HW0}\beg{split} \mu(f^2\log f^2) &\le  r\mu(|\nn f|^2) +\ff{(Kr+2)W_2^\rr(f^2\mu,\mu)}{2r} \Big(W_2^\rr(f^2\mu,\mu)\land \ff{\ss {rn}}{2\ss 2}\Big) \\
&\quad +\ff{\ss n (Kr+2)} {4\ss{2r}} \Big(  W_2^\rr(f^2\mu,\mu) - \ff{\ss {rn}}{2\ss 2} \Big)^+,\ \ r\in (0,\infty)\cap \Big(0,\ff 2 {K^-}\Big],\end{split}\end{equation} where $K^-:= \max\{0, -K\}$. Consequently,
\beq\label{HWI} \beg{split} \mu(f^2\log f^2)\le  & 2W_2^\rr(f^2\mu,\mu)\ss{\mu(|\nn f|^2)} +\ff K 2 W_2^\rr(f^2\mu,\mu)^2\\
  -&\ff{KW_2^\rr(f^2\mu,\mu)+2\ss{\mu(|\nn f|^2)}}{2\ss{W_2^\rr(f^2\mu,\mu)}}\Big( \ss{W_2^\rr(f^2\mu,\mu)} - \ff{ \ss{n}}{2\ss 2 \mu(|\nn f|^2)^{1/4}}\Big)^+.\end{split}\end{equation}\end{cor}

It was proved in \cite{OV} and \cite{BGL} that (\ref{C}) (i.e. (\ref{CD}) for $n=\infty$) implies

$$\mu(f^2\log f^2)\le  2W_2^\rr(f^2\mu,\mu)\ss{\mu(|\nn f|^2)} +\ff K 2 W_2^\rr(f^2\mu,\mu)^2 $$   for all $f\in C^1(M)$ with $\mu(f^2)=1.$ According to
(\ref{HWI}), the dimension $n$ contributes to a negative term in the right-hand side since $KW_2^\rr(f^2\mu,\mu)+2\ss{\mu(|\nn f|^2)}\ge 0$ as explained in the proof of (\ref{HWI}). But this inequality is incomparable with the Sobolev type WHI inequality derived in \cite{W08}.

\

Next, we  consider the first non-trivial eigenvalue (i.e. the spectral gap) of $L$. To this end, let $Z=\nn V$ for
some $V\in C^2(M)$ such that

$$\mu(\d x):=\e^{V(x)}\d x$$ is a probability measure, where $\d x$ stands for the Riemannian volume measure on the manifold.
In this case the Friedrich extension of $(L,\D_0)$ gives rise to a negatively definite self-adjoint operator on $L^2(\mu)$, whose spectral gap   can be characterized  as

$$\ll_1=\inf\{\mu(|\nn f|^2): f\in C_0^1(M), \mu(f)=0, \mu(f^2)=1\}.$$
The following lower bound of $\ll_1$ is a simple consequence of Theorem \ref{T1.1} (2). This estimate is well known as the Lichnerowicz estimate \cite{Lich} for $Z=0$, and was extended to $Z\ne 0$ by Bakry and Qian \cite{BQ00}.

\beg{cor}[\cite{Lich, BQ00}]\label{C1.2} Let $Z=\nn V$ such that $\mu(\d x):=\e^{V(x)}\d x$ is a probability measure. If $(\ref{CD})$ holds
for some $K<0$ and $n>1$, then

$$\ll_1\ge \ff{n(-K)}{n-1}.$$\end{cor}

\

Finally, we consider the transportation inequality of $P_t$ deduced from (\ref{CD}).  According to \cite{RS}, (\ref{C}) implies

\beq\label{CT} W_p^\rr(\mu_1 P_t,\mu_2 P_t)\le \e^{Kt}W_p^\rr(\mu_1,\mu_2),\ \ t\ge 0, \mu_1,\mu_2\in \scr P(M) \end{equation} for any $p\ge 1$. Using (\ref{CD}) we prove the following inequalities (\ref{CT''}) and (\ref{CT'}). Comparing with (\ref{CT}), when $p=1$ (\ref{CT'})   has better long time behavior for $K<0$ while (\ref{CT''}) is stronger for $K>0$. In fact,  since

$$H(r):= \ff 2 {\ss{K/(n-1)}}\sinh\Big[\ff r 2 \ss{K/(n-1)}\Big],\ \ r\ge 0$$ is convex with $H'(r)>1$ for $r>1$, due to the Jensen inequality,  (\ref{CT''})   implies that

\beg{equation*}\beg{split} W_1^\rr(\mu_1 P_t,\mu_2 P_t) &\le H^{-1}  \big( W_1^{\tt\rr}(\mu_1 P_t, \mu_2 P_t) \big)
\le H^{-1} \big(\e^{ Kt} H(W_1^{\tt\rr}(\mu_1,\mu_2))\big)\\
&<\e^{Kt} H^{-1}\circ H(W_1^{\tt\rr}(\mu_1,\mu_2))=\e^{Kt}W_1^{\tt\rr}(\mu_1,\mu_2),\ \ \ t>0, \mu_1\ne \mu_2.\end{split}\end{equation*}

\beg{prp}\label{PP} Assume that $(\ref{CD})$ holds and let
$$\tt \rr(x,y)=  \beg{cases} \ff{2}{\ss{-K/(n-1)}} \sin\Big[\ff {\rr(x,y)} 2 \ss{-K/(n-1)}\Big], &{\rm if}\ K<0,\\
\rr(x,y), &{\rm if}\ K=0,\\
\ff{2}{\ss{K/(n-1)}} \sinh\Big[\ff {\rr(x,y)} 2 \ss{K/(n-1)}\Big], &{\rm if}\ K>0.\end{cases}$$ Then for any $p\ge 1$,
\beq\label{CT''} W_p^{\tt\rr} (\mu_1 P_t, \mu_2 P_t)\le \e^{Kt} W_p^{\tt\rr}(\mu_1,\mu_2),\ \ t\ge 0, \mu_1,\mu_2\in\scr P(M).\end{equation}
If $K<0$ then \beq\label{CT'} W_1^{\tt\rr} (\mu_1 P_t, \mu_2 P_t)\le   \exp\Big[\ff{n K}{n-1} t\Big] W_1^{\tt\rr}(\mu_1,\mu_2),  \ \ t\ge 0, \mu_1,\mu_2\in\scr P(M).\end{equation}  \end{prp}

\section{Proofs} According to the proof of \cite[Theorem 1.1]{W10}, the reflection at a convex boundary does not make any trouble for our proofs. So, for simplicity, we shall only consider the case without boundary. In this case, the the proofs of (1)-(5) in Theorem \ref{T1.1} are more or less standard according to the semigroup argument of Bakry, Emery and Ledoux. Our proof of equivalence between (6) and (\ref{CD}) is however highly technical.

 \beg{proof}[Proof of Theorem \ref{T1.1}] By the Jensen inequality, (2) follows from (1) immediately. So, it suffices to show that (\ref{CD}) implies (1), (2) implies (3) and (4), each of (3) and (4) implies (\ref{CD}),  (5) is equivalent to (\ref{CD'}), (2) implies (6), and (6) implies (\ref{CD}). Below we prove these implications respectively.

\paragraph{(\ref{CD}) implies (1).} By (\ref{CD}) we have

\beg{equation*}\beg{split} \ff{\d}{\d s} P_s |\nn P_{t-s} f|^2 &= P_s\big\{L|\nn P_{t-s} f|^2 -2\<\nn P_{t-s}f, \nn L P_{t-s} f\>\big\} \\
&\ge -2K P_s|\nn P_{t-s}f|^2 + \ff 2 n P_s (P_{t-s}Lf)^2,\ \ s\in [0,t].\end{split}\end{equation*} By the Gronwall lemma, this implies (1) immediately. \paragraph{(2) implies (3) and (4).} Obviously, we have

\beq\label{2.1} \ff{\d }{\d s} P_s (P_{t-s} f)^2 =P_s\{L(P_{t-s}f)^2 -2 (P_{t-s}f)  LP_{t-s}f\}=2P_s|\nn P_{t-s}f|^2.\end{equation} Next, according to (2) and noting that $P_s(P_{t-s}Lf)^2\ge (P_t Lf)^2,$ we have

\beg{equation*}\beg{split} & P_s|\nn P_{t-s}f|^2 \le \e^{2K(t-s)} P_t|\nn f|^2 -\ff{\e^{2K(t-s)}-1}{Kn} P_s(P_{t-s} Lf)^2,\\
& P_s|\nn P_{t-s}f|^2 \ge \e^{-2Ks} |\nn P_t  f|^2 + \ff{1-\e^{-2Ks}}{Kn}  (P_{t} Lf)^2.\end{split}\end{equation*} Combining these with (\ref{2.1}) respectively and integrating w.r.t. $\d s$ over $[0,t]$, we prove (3) and (4).
\paragraph{  (3) or (4) implies (\ref{CD}).}  For small $t>0$ we have

\beg{equation*}\beg{split} &P_t f^2 = f^2 t L f^2 +\ff {t^2} 2 L^2 f^2 + \circ(t^2),\\
&(P_t f)^2 = \Big(f+ t Lf +\ff {t^2} 2 L^2 f+\circ(t^2)\Big)^2 = f^2 +t^2(Lf)^2 + 2t f Lf +t^2 fL^2 f+\circ(t^2).\end{split}\end{equation*} So,

\beq\label{2.2} P_t f^2 -(P_tf)^2 = 2 t |\nn f|^2 +t^2 \{2\<\nn L f, \nn f\> +L|\nn f|^2\} +\circ(t^2).\end{equation} On the other hand,

\beg{equation*}\beg{split} \ff{\e^{2Kt}-1} K P_t |\nn f|^2 &= \{2t + 2Kt^2 +\circ(t^2)\}\cdot\{|\nn f|^2 +t L|\nn f|^2 +\circ(t)\}\\
&= 2t |\nn f|^2 + 2t^2 \{L|\nn f|^2 +K |\nn f|^2\}+\circ(t^2).\end{split}\end{equation*} Moreover, it is easy to see that

$$\ff{\e^{2Kt}-2Kt-1}{K^2 n} (P_t Lf)^2 = \ff 2 n t^2 (Lf)^2 +\circ(t^2).$$ Combining these with (\ref{2.2}), we see that (3) implies

$$2t^2 \Big\{\ff 1 2 L|\nn f|^2 -\<\nn Lf, \nn f\> + K|\nn f|^2 +\ff{(Lf)^2} n\Big\}+\circ(t^2) \ge 0.$$ Therefore, (\ref{CD}) holds.

Next, it is easy to see that

\beg{equation*}\beg{split} &\ff{1-\e^{-2Kt}} K |\nn P_t|^2 +\ff{\e^{-2Kt}-1+2Kt}{K^2n} (P_tLf)^2\\
&= \{2t -2Kt^2 +\circ(t^2)\}\cdot |\nn f +t \nn Lf +\circ(t)|^2 + \ff{2t^2} n (Lf)^2 +\circ(t^2)\\
&= 2t|\nn f|^2 + 2t^2\Big\{2   \<\nn f, \nn Lf\> + \ff {(Lf)^2} n - K|\nn f|^2\Big\}+\circ (t^2).\end{split}\end{equation*}
Combining this with (\ref{2.2}) and (4) we prove (\ref{CD}).
\paragraph{(5) is equivalent to (\ref{CD'}).}  Using $\ss{|\nn P_{t-s}f|^2+\vv}$ to replace $|\nn P_{t-s}f|$ and letting $\vv\to 0$, in the following calculations we may assume that $|\nn P_{t-s}f|$ is positive and smooth, so that

\beg{equation}\label{AA}\beg{split} \ff{\d}{\d s} P_s|\nn P_{t-s}f| &= P_s\bigg\{L|\nn P_{t-s}f|- \ff{\<\nn L P_{t-s} f, \nn P_{t-s}f\>}{|\nn P_{t-s}f|}\bigg\}\\
&= P_s\bigg\{\ff{\ff 1 2 L|\nn P_{t-s}f|^2 -\<\nn LP_{t-s}f, \nn P_{t-s} f\> -\big|\nn |\nn P_{t-s}f|\big|^2}
{|\nn P_{t-s}f|}\bigg\}.\end{split}\end{equation} Since
\beq\label{HH} \beg{split}&\ff 1 2L|\nn f|^2 -\<\nn Lf,\nn f\> =\Ric(\nn f,\nn f)-\<\nn_{\nn f} Z, \nn f\> +\|{\rm Hess}_f\|_{HS}^2,\\
&\big|\nn |\nn f|\big|^2 = \bigg\|{\rm Hess}_f\Big(\ff{\nn f}{|\nn f|},\cdot\Big)\bigg\|^2 \le \|{\rm Hess}_f\|_{HS}^2,\end{split}\end{equation} it follows from (\ref{CD'}) and (\ref{AA}) that

$$\ff{\d}{\d s}P_s|\nn P_{t-s}f| \ge -K P_s|\nn P_{t-s}f| +\ff 1 {n-d} P_s \ff{\<Z, \nn P_{t-s}f\>^2}{|\nn P_{t-s}f|}.$$ This implies (5).

On the other hand, since when $t=0$ the equality in (5) holds, one may take derivatives at $t=0$ for both sides of (5) to derive at points  such that $|\nn f|>0$

$$ K|\nn f|+L|\nn f| \ge \ff{\<\nn L f, \nn f\>}{|\nn f| } +\ff {\<Z,\nn f\>^2} {(n-d)|\nn f|}.$$ This implies

$$\ff 1 2 L|\nn f|^2 -\<\nn Lf,\nn f\>\ge -K|\nn f|^2 +\ff{\<Z,\nn f\>^2}{n-d}.$$ Combining this with (\ref{HH}) we obtain

$$\Ric(\nn f,\nn f)-\<\nn_{\nn f}Z,\nn f\>\ge -K |\nn f|^2 + \ff{\<Z,\nn f\>^2}{n-d},\ \ f\in C^\infty(M),$$ which is equivalent to (\ref{CD'}).
\paragraph{(2) implies (6).} By the monotone class theorem, we may assume that $f\in C^2(M)$ which is constant outside a compact set. Let $\gg: [0,1]\to M$ be the minimal geodesic from $x$ to $y$, and let

$$h(s)= \ff{\int_0^s \e^{-2K\varphi(r)}\d r}{\int_0^t \e^{-2K\varphi(r)}\d r},\ \ \ s\in [0,t].$$ By (2) we have

\beg{equation*}\beg{split} &\ff{\d}{\d s}P_{\varphi(s)}\log P_{t-s} f(\gg_{h(s)})\\
 &= P_{\varphi(s)}\Big\{\varphi'(s)L\log P_{t-s} f- \ff{LP_{t-s} f}{P_{t-s} f}\Big\}(\gg_{h(s)}) + h'(s)\<\dot \gg_{h(s)}, \nn P_{\varphi(s)}\log P_{t-s} f(\gg_{h(s)})\>\\
&\le P_{\varphi(s)} \big\{(\varphi'(s)-1) L\log P_{t-s}f- |\nn \log P_{t-s}f|^2\}(\gg_{h(s)}) \\
&\quad+\big\{|h'(s)|\rr(x,y)\big\}|\nn P_{\varphi(s)} \log P_{t-s} f|(\gg_{h(s)})\\
&\le  \Big\{|\varphi'(s)-1|\cdot\big|P_{\varphi(s)}L\log P_{t-s}f\big|  -
\ff{1-\e^{-2K\varphi(s)}}{Kn} (P_{\varphi(s)}L\log P_{t-s} f)^2\Big\}(\gg_{h(s)})\\
&\quad+\Big\{\rr(x,y)h'(s)\big|\nn P_{\varphi(s)}\log P_{t-s}f\big|- \e^{-2K\varphi(s)}
 \big|\nn P_{\varphi(s)}\log P_{t-s}f\big|^2\Big\}(\gg_{h(s)})\\
&\le \ff{\e^{2K\varphi(s)}\rr(x,y)^2 h'(s)^2}{4}+ \ff{Kn(\varphi'(s)-1)^2}{4(1-\e^{-2K\varphi(s)})}.\end{split}\end{equation*}This completes the proof by integrating w.r.t. $\d s$ over $[0,t]$.
\paragraph{(6) implies (\ref{CD}).} For fixed $x\in M$ and strictly positive $f\in C^\infty(M)$ which is constant outside a compact set. Let

$$\varphi(s) = s+ \ff {2L(\log f)(x)}{n} s^2,\ \ \gg_s= \exp[-2s \nn\log f(x)],\ \ s\ge 0.$$ According to (6), for small $t>0$ we have
\beq\label{HT} P_{\varphi(t)}(\log f)(x) \le \log P_t f(\gg_t)+\ff{t^2|\nn \log f|^2(x)}{\int_0^t \e^{-2K\varphi(s)}\d s}
+\ff {Kn} 4 \int_0^t \ff{(\varphi'(s)-1)^2}{1-\e^{-2K\varphi(s)}}\d s.\end{equation}
According to (3.3) in \cite{W10} and noting that $\varphi(t)^2 = t^2 +\circ(t^2)$, we have
 \beq\label{6-1} \beg{split} & P_{\varphi(t)} (\log f)(x)= \log f(x) + \varphi(t) L\log f(x) + \circ(t^2)\\
&+ \ff {t^2} 2\Big\{\ff{L^2 f}f -\ff{(Lf)^2}{f^2} -\ff {2  \<\nn Lf,\nn f\>}{f^2} -\ff{L|\nn f|^2}{f^2}+\ff{4|\nn f |^2Lf} {f^3}-\ff{6|\nn f|^4}{f^4}\Big\}.\end{split}\end{equation}
Moreover, according to   line 10 on page 310 in \cite{W10} and noting that we do not assume $\Hess_f(x)=0$,
\beq\label{6-2} \beg{split} &\log P_t  f (\gg_t)= \log f(x) + t\big\{ L\log f(x) -|\nn\log f|^2\big\}(x) + \circ(t^2)\\
&+ \ff {t^2} 2\Big\{\ff{L^2 f}f -\ff{(Lf)^2}{f^2} -\ff {4  \<\nn Lf,\nn f\>}{f^2} +\ff{4  |\nn f|^2 Lf} {f^3}   - \ff{4|\nn f|^4}{f^4}+\ff{4\Hess_f(\nn f,\nn f)}{f^3}\Big\}.\end{split}\end{equation} Finally, since it is easy to see that
$$\lim_{t\to 0} \ff{Kn}{4t^2}\int_0^t \ff{(\varphi'(s)-1)^2}{1-\e^{-2K\varphi(s)}}\d s = (L\log f)^2(x),$$ we have
\beq\label{6-3}\ff{Kn}{4}\int_0^t \ff{(\varphi'(s)-1)^2}{1-\e^{-2K\varphi(s)}}\d s= t^2(L\log f)^2(x) +\circ(t^2).\end{equation}
Substituting (\ref{6-1}), (\ref{6-2}) and (\ref{6-3}) into (\ref{HT}), and noting that
$$(\varphi(t)-t)L(\log f)(x)= \ff{2t^2}n (L\log f)^2(x),$$
we arrive at
\beg{equation*}\beg{split}&\ff 1 t \bigg(1-\ff t {\int_0^t \e^{-2K\varphi(s)}\d s}\bigg)|\nn\log f|^2(x) +\ff{(L\log f)^2 (x)}{n}\\
 &\le \ff 1  2 \bigg(\ff{L|\nn f|^2-2\<\nn L f,\nn f\>} {f^2}+\ff{2|\nn f|^4}{f^4} +\ff{4|\Hess_f(\nn f,\nn f)|}{f^3}\bigg)(x)+\circ(1).\end{split}\end{equation*} Letting $t\to 0$ and multiplying both sides by $f^2$, we obtain
\beg{equation*}\beg{split}&-K|\nn f|^2(x) + \ff {(Lf - |\nn f|^2/f)^2(x)}{n}\\
 &\le \bigg(\ff 1 2  L|\nn f|^2- \<\nn L f,\nn f\>  +\ff{ |\nn f|^4}{f^2} +\ff{2|\Hess_f(\nn f,\nn f)|}{f }\bigg)(x).\end{split}\end{equation*}
Replacing $f$ by $f+m$ and letting $m\to\infty$, this implies that
$$-K|\nn f|^2(x) +\ff{(Lf)^2(x)}n \le \ff 1 2  L|\nn f|^2(x)- \<\nn L f,\nn f\> (x).$$ Therefore, (\ref{CD}) holds.
\end{proof}

\beg{proof}[Proof of Corollary \ref{C1.3}]  Let $t_0\in (0,t)$. Taking

$$\varphi(r) =r\land \ff t 2  + \ff{t+2s}{t} \Big(r-\ff t2\Big)^+,\ \ r\in [0,t],$$ we have
\beg{equation*}\beg{split} \int_0^t \e^{-2K\varphi(r)}\d r&= \ff{1-\e^{-Kt}}{2K}+\ff{t(\e^{-Kt}-\e^{-2K(t+s)})}{2K(t+s)}\\
&\ge \ff{t(1-\e^{-2K(t+s)})}{2K(t+2s)},\end{split}\end{equation*} and
$$ K\int_0^t \ff{(\varphi'(r)-1)^2}{1-\e^{-2K\varphi(r)}}\,\d r= \ff{4Ks^2}{t^2} \int_{t/2}^t \ff{\d r}{1-\exp[-\ff{2K(t+2s)} t(r-\ff t 2)-Kt]}
 \le \ff{2Ks^2}{t(1-\e^{-Kt})}.$$ Thus, (\ref{H1}) follows from  (6).

Next, applying Lemma (6) for $t+s$ in place of $t$ and taking $\varphi(r)= r\land t$, we prove (\ref{H2}).
\end{proof}

\beg{proof}[Proof of Corollary \ref{C1.3'}]  When $s=0$,   (\ref{H1}) and (\ref{H2}) hold for $n=\infty$ (see \cite{W10}). Applying e.g. (\ref{H1}) to $s=0$ and $f(z):= p_{t}(y,z)\land m+\vv$ for $m,\vv>0$ and letting $m\to\infty, \vv\to 0$, we obtain

$$\int_Mp_{t}(y,z)\log p_{t}(y,z)\mu(\d z)\le \log p_{2t}(x,y) +\ff{K\rr(x,y)^2}{2(1-\e^{-2K(t)})}.
$$ Since $\mu$ is a probability measure and $\int_Mp_{t+s}(y,z)\mu(\d z)=1,$ by the Jensen inequality this implies

$$p_{2t}(x,y)\ge  \exp\bigg[-\ff{K\rr(x,y)^2}{2(1-\e^{-2Kt})}\bigg].$$ Replacing $t$ by $\ff t2$, we prove the desired heat kernel lower bound. \end{proof}

\beg{proof}[Proof of Corollary \ref{C1.4}] Applying (\ref{H2}) for $P_tf^2+\vv$ in place of $f$ and letting $\vv\to 0$, we obtain

$$ (P_t\log P_t f^2)(x)  \le \log P_{2t+s}f^2(y) +\ff{\rr(x,y)^2 K}{2(1-\e^{-2Kt})+4s K\e^{-2Kt}} +\ff{Kns}{4(1-\e^{-2Kt})},\ \ s\ge 0.$$ Let $\pi\in\scr C(f^2\mu,\mu)$ be the optimal coupling for $W_2^\rr(f^2\mu,\mu)$, integrating both sides w.r.t. $\pi$ and noting that due to the Jensen inequality and
$\mu(f^2)=1$ it follows that $\mu(\log P_{2t+s}f^2)\le 0$, we arrive at

\beq\label{2.2} \mu((P_tf^2)\log P_t f^2) \le \ff {W_2^2K}{2(1-\e^{-2Kt})+ 4 s \e^{-2Kt}K}+\ff{Kns}{4(1-\e^{-2Kt})},\end{equation} where and in the remainder of the proof, $W_2$ stands for $W_2^\rr(f^2\mu,\mu)$ for simplicity. On the other hand, it is well known that (\ref{CD}) (indeed, (\ref{C}) implies

$$P_tf^2\log f^2\le (P_t f^2)\log P_t f^2 +\ff{\e^{2Kt}-1}K P_t |\nn f|^2.$$ Integrating both sides w.r.t. $\mu$ and using (\ref{2.2}) we obtain

$$\mu(f^2\log f^2)\le \ff{\e^{2Kt}-1}K \mu(|\nn f|^2)+ \ff {W_2^2K}{2(1-\e^{-2Kt})+ 4 s \e^{-2Kt}K}+\ff{Kns}{4(1-\e^{-2Kt})}.$$
Letting $r= 2(\e^{2Kt}-1)/K$ which runs over all $(0,\ff 2 {K^-})$ as $t$ varies in $(0,\infty)$, and using $rs$ to replace $s$, we get
$$\mu(f^2\log f^2)\le r \mu(|\nn f|^2) +(Kr+2)\Big\{\ff{W_2^2}{2(1+4 s)r}+\ff{ns}4\Big\}, \ \ 0<r\le \ff 2 {K^-}, s>0.$$ Taking

$$s=\ff 1 {4} \Big(\ff{2\ss 2 W_2}{\ss{rn}} -1\Big)^+,$$ we prove (\ref{HW0}).

To prove (\ref{HWI}), let
$$ \dd = \mu(|\nn f|^2),\ \ \ r= \ff{W_2}{\ss\dd}.$$   Since  according to \cite{BE, OV} one has
$$\ff {K^-} 2 W_2^\rr (f^2\mu,\mu)^2\le \mu(f^2\log f^2)\le \ff 2 {K^-}\mu(|\nn f|^2),$$ it is clear that $r\le \ff 2{K^-}$. Thus,   (\ref{HW0}) applies to this specific $r$. Therefore, (\ref{HWI}) follows by noting that
\beg{equation*}\beg{split}
&r\dd+ \ff{(Kr+2)W_2}{2r} \Big( W_2 \land \ff{\ss {rn}}{2\ss 2}\Big) +\ff{\ss n (Kr+2)} {4\ss{2r}} \Big( W_2 -\ff{\ss {rn}}{2\ss 2}\Big)^+\\
&= \dd r +\ff{(Kr+2)W_2^2}{2r} - \ff{(Kr+2)W_2}{2r}\Big( W_2 -\ff{\ss {rn}}{2\ss 2}\Big)^+ +\ff{\ss n (Kr+2)} {4\ss{2r}} \Big( W_2 -\ff{\ss {rn}}{2\ss 2}\Big)^+\\
&= \dd r +\Big(\ff K 2 + \ff 1 r \Big)W_2^2 - \ff{Kr+2}{2r} {\Big( W_2 -\ff{\ss {rn}}{2\ss 2}\Big)^+}^2\\
&= 2W_2\ss\dd +\ff K 2 W_2^2 -\ff{KW_2 +2\ss\dd}{2\ss{W_2}}{\Big(\ss {W_2} -\ff{\ss n}{2\ss 2 \dd^{1/4}}\Big)^+}^2.\end{split}\end{equation*} \end{proof}

\beg{proof}[Proof of Corollary \ref{C1.2}] Since $K<0$, the manifold is compact (cf. \cite{L}). In this case the spectrum of $L$ is discrete so that $\ll_1>0$ and there exists an eigenfunction $f$ with $\mu(f^2)=1$ and $Lf=-\ll_1f$. By Theorem \ref{T1.1}(2) we have

$$\mu(|\nn P_t f|^2)\le \e^{2Kt}\mu(|\nn f|^2) -\ff{\e^{2Kt}-1}{Kn} \mu((P_t Lf)^2),\ \ t>0.$$ For  $f$ being the above mentioned eigenfunction, this implies

$$\ll_1\e^{-2\ll_1 t}\le \ll_1 \e^{2Kt} -\ll_1^2 \e^{-2\ll_1 t}\ff{\e^{2Kt}-1}{Kn},\ \ t>0.$$ Equivalently,

$$ \ff{\e^{2(K+\ll_1)t}-1} t \ge \ll_1  \ff{\e^{2Kt}-1}{Knt},\ \ t>0.$$ Letting $t\to 0$ we obtain the desired lower bound of $\ll_1.$\end{proof}

\beg{proof}[Proof of Proposition \ref{PP}] Since the assertion for $K=0$ follows from that for $K>0$ by letting $K\to 0$, below we only prove the desired inequality for  $K<0$ and $K>0$ respectively.

(a) Let $K<0$. Take $\pi\in \C(\mu_1,\mu_2)$ such that $W_1^{\tt\rr}(\mu_1,\mu_2)= \pi(\tt\rr)$, and let $(X_0,Y_0)$ be an $M\times M$-valued random variable with distribution $\pi$. Let $(X_t,Y_t)$ be the coupling by reflection of the $L$-diffusion process with initial data $(X_0,Y_0)$. This coupling was initiated by Kendall \cite{K} and Cranston \cite{Cr} (see
\cite[\S 2.1]{W05} for a complete construction). We have (see \cite{CW97} or \cite[Theorem 2.1.1]{W05})
\beq\label{D0} \d\rr(X_t,Y_t)\le 2\ss 2\,\d b_t +I_Z(X_t, Y_t)\d t\end{equation}
for a one-dimensional Brownian motion $b_t$ and
\beq\label{0}I_Z(x,y):= I(x,y)+ \<Z,\nn\rr(\cdot, y)\>(x)+ \<Z,\nn \rr(x,\cdot)\>(y),\end{equation} where letting $\gg: [0,\rr(x,y)]\to M$ be the minimal geodesic from $x$ to $y$ and $\{J_i\}_{i=1}^{d-1}$ the Jacobi fields along $\gg$ such that at points $x,y$ they together with $\dot\gg$ consist of an orthonormal  basis of the tangent space, we have
$$I(x,y)= \sum_{i=1}^{d-1} \int_0^{\rr(x,y)} \big(|\nn_{\dot\gg} J_i|^2 -\<\scr R(\dot\gg, J_i)\dot\gg, J_i\>\big)_s\d s,$$ where $\scr R$ is the   curvature tensor on $M$.

To calculate $I(x,y)$, let us fix points $x\ne y$ and simply denote $\rr=\rr(x,y)$. Let $\{U_i\}_{i=1}^{d-1}$ be constant vector fields along $\gg$ such that $\{\dot \gg, U_i: 1\le i\le d-1\}$ is an orthonormal basis. By the index lemma, for any $f\in C^1([0,\rr])$ with $f(0)=f(\rr)=1$, we have
\beq\label{D1}\beg{split}  I(x,y)&\le \sum_{i=1}^{d-1} \int_0^\rr \big(|\nn_{\dot\gg} f U_i|^2 -f^2\<\scr R(U_i, \dot\gg) \dot\gg, U_i\>\big)_s\d s\\
&= \int_0^\rr \big\{(d-1) f'(s)^2 -f(s)^2 \Ric(\dot\gg,\dot\gg)_s\big\}\d s.\end{split}\end{equation} On the other hand, since $f(0)=f(\rr)=1$,
\beg{equation*}\beg{split} &\<Z,\nn\rr(\cdot, y)\>(x)+ \<Z,\nn \rr(x,\cdot)\>(y)= \int_0^\rr\ff{\d}{\d s} \big\{f(s)^2\<\dot\gg, Z\circ\gg\>_s\big\}\d s\\
&=\int_0^\rr \big\{2(ff')(s) \<\dot\gg, Z\circ\gg\>_s +f(s)^2 \<\nn_{\dot\gg}Z\circ\gg,\dot\gg\>_s\big\}\d s\\
&\le \int_0^\rr\bigg\{\ff{f(s)^2\<\dot\gg, Z\circ\gg\>_s^2}{n-d} +(n-d)f'(s)^2 + f(s)^2 \<\nn_{\dot\gg}Z\circ\gg,\dot\gg\>_s\bigg\}\d s.\end{split}\end{equation*} Combining this with (\ref{D1}), (\ref{0})  and  (\ref{CD'}), we obtain
\beq\label{D2} I_Z(x,y)\le \int_0^\rr\big[(n-1) f'(s)^2 +Kf(s)^2\big]\d s.\end{equation}
Taking
$$f(s)= \tan\Big(\ff{\rr} 2\ss{-K/(n-1)}\sin\Big(\ss{-K/(n-1)}\,s\Big)+ \cos\Big(\ss{-K/(n-1)}\,s\Big)$$ for $s\in [0,\rr],$  we obtain
\beq\label{DW} I_Z(x,y)\le -2 \ss{-K(n-1)}\tan\Big(\ff{\rr} 2 \ss{-K/(n-1)}\Big).\end{equation} Therefore, it follows from (\ref{D0}) and the It\^o formula  that
$$\d\tt\rr(X_t,Y_t)\le \d M_t+ \ff{nK}{n-1} \tt\rr(X_t,Y_t)\d t$$ holds for some martingale $M_t$. Thus,
$$W_1^{\tt\rr}(\mu_1 P_t, \mu_2 P_t)\le \E\tt\rr(X_t,Y_t) \le\exp\Big[\ff{n K}{n-1} t\Big] \E \tt\rr(X_0, Y_0)=
\exp\Big[\ff{n K}{n-1} t\Big] W_1^{\tt\rr}(\mu_1,\mu_2).$$

(b) When $K>0$, we take
\beg{equation*}\beg{split} f(s) = &\cosh\Big(\ff{\rr} 2\ss{K/(n-1)}\sinh\Big(\ss{K/(n-1)}\,s\Big)\\
&+ \ff{1-\cosh(\rr\ss{K/(n-1)})}{\sinh(\rr\ss{K/(n-1)})}\sinh \big(s\ss{K/(n-1)}\big),\ \ s\in [0,\rr].
 \end{split}\end{equation*} It follows from (\ref{D2}) that
 $$I_Z(x,y) \le 2\ss{K(n-1)}\tanh\Big(\ff{\rr(x,y)} 2 \ss{K/(n-1)}\Big).$$ Combining this with (\ref{DW}), we obtain
\beq\label{D3}I_Z(x,y)=\beg{cases} 2\ss{K(n-1)}\tanh\Big(\ff{\rr(x,y)} 2 \ss{K/(n-1)}\Big), &\text{if} \ K>0;\\
-2 \ss{-K(n-1)}\tan\Big(\ff{\rr(x,y)} 2 \ss{-K/(n-1)}\Big), &\text{if}\ K<0.\end{cases}\end{equation} Now, let $(X_0,Y_0)$ have distribution $\pi$ such that $\pi(\tt\rr^p)= W_p^{\tt\rr}(\mu_1,\mu_2)^p.$ Using the coupling by parallel displacement rather than by reflection, we have
(see \cite[Proof of Proposition 2.5.1]{W05} or \cite{ATW})
$$\d\rr(X_t,Y_t)\le I_Z(X_t,Y_t)\d t.$$ Combining this with   (\ref{D3}) we conclude that
$$\d\tt \rr(X_t,Y_t) \le \e^{K t} \tt\rr(X_t,Y_t).$$ Therefore,
$$W_p^{\tt\rr}(\mu_1P_t,\mu_2P_t)\le (\E \tt\rr(X_t,Y_t)^p)^{1/p} \le \e^{Kt} (\E \tt\rr(X_0,Y_0)^p)^{1/p}=\e^{Kt} W_p^{\tt\rr}(\mu_1,\mu_2).$$\end{proof}

\beg{thebibliography}{99}

\bibitem{ATW} M. Arnaudon,  A. Thalmaier,    F.-Y. Wang,
 \emph{Harnack inequality and heat kernel estimates on
manifolds with curvature unbounded below,}  Bull. Sci. Math.
130(2006), 223--233.

\bibitem{B} D. Bakry, \emph{On Sobolev and logarithmic
 Sobolev inequalities for
Markov semigroups,} New Trends in Stochastic Analysis, 43--75,
 World Scientific, 1997.

 \bibitem{BE}   D. Bakry,  M. Emery, \emph{Hypercontractivit\'e de
semi-groupes de diffusion}, C. R. Acad. Sci. Paris. S\'er. I Math.
299(1984), 775--778.

\bibitem{BQ99} D. Bakry, Z. Qian, \emph{Harnack inequalities on a manifold with positive or negative Ricci curvature,}
 Rev. Mat. Iberoamericana 15 (1999),  143--179.

 \bibitem{BQ00} D. Bakry, Z. Qian, \emph{Some new results on eigenvectors via dimension, diameter and Ricci curvature,}
 Adv. Math. 155(2000), 98--153.

\bibitem{BGL} S. G. Bobkov, I.  Gentil and M.  Ledoux,
\emph{Hypercontractivity of Hamilton-Jacobi equations,} J. Math.
Pures Appl. 80(2001), 669--696.

\bibitem{CW97}  M.-F. Chen, F.-Y. Wang, \emph{Application of coupling method to the first eigenvalue on manifold,} Sci. in China (A) 40(1997), 384--394.

\bibitem{Cr} M. Cranston, \emph{Gradient estimates on manifolds using coupling}, J. Funct.
  Anal. 99 (1991), 110--124.

\bibitem{K} W.\,S. Kendall, \emph{Nonnegative {R}icci curvature and the {B}rownian
  coupling property}, Stochastics  19 (1986),   111--129.

\bibitem{L} M. Ledoux, \emph{The geometry of Markov diffusion generators,} Ann. Facu. Sci. Toulouse
9(2000), 305--366.

\bibitem{LY} P. Li, S.-T. Yau, \emph{On the parabolic kernel of the Schr\"odinger operator,} Acta Math.
156 (1986), 153--201.

\bibitem{Lich} A. Lichnerowicz, \emph{G\'eom\'etrie des Groupes des Transformations,} Paris: Dnnod, 1958.

\bibitem{OV} F.  Otto,  C. Villani, \emph{Generalization of an inequality
by Talagrand and links with the logarithmic Sobolev inequality,} J.
Funct. Anal.  173(2000), 361--400.

\bibitem{Q} Z. Qian, \emph{A comparison theorem for an elliptic operator,}
 Potential Anal. 8 (1998),  137--142.

\bibitem{RS} M.-K. von Renesse,  K.-T.  Sturm, \emph{Transport inequalities, gradient estimates, entropy,
and Ricci curvature,}  Comm. Pure Appl. Math. 58 (2005),  923--940.

\bibitem{RW}  M. R\"ockner, F.-Y. Wang, \emph{Log-Harnack  inequality for stochastic differential equations in Hilbert spaces and its consequences, } Infin. Dimens. Anal. Quant. Probab.  Relat. Topics 13(2010), 27--37.

\bibitem{S} L. Saloff-Coste, \emph{Convergence to equilibrium and elliptic operator,} Potential Anal. 8(1998), 137--142.

\bibitem{St} K.-T. Sturm, \emph{Heat kernel bounds on manifolds,} Math. Ann. 292(1992), 149--162.

\bibitem{W04} F.-Y. Wang, \emph{Equivalence of dimension-free Harnack inequality and curvature condition,}
 Integral Equations Operator Theory 48 (2004),  547--552.

 \bibitem{W05} F.-Y. Wang, \emph{Functional Inequalities, Markov Semigroups and Spectral Theory,} Science Press, Beijing, 2005.

 \bibitem{W08}	F.-Y. Wang, \emph{Generalized transportation-cost inequalities and applications,} Potential Anal. 28(2008), 321--334.

\bibitem{W10}	F.-Y. Wang, \emph{Harnack inequalities on manifolds with boundary and applications,}  J. Math. Pures Appl. 94(2010), 304--321.

\end{thebibliography}

\end{document}